\documentclass{ws-ijbc}
\usepackage{}
\usepackage{mathrsfs}
\usepackage{multicol}
\usepackage{ws-rotating}     
\begin{document}

\catchline{}{}{}{}{} 

\markboth{ANHUI GU}{Random Attractors of Stochastic Lattice
Dynamical Systems Driven by Fractional Brownian Motions}

\title{RANDOM ATTRACTORS OF STOCHASTIC LATTICE DYNAMICAL
SYSTEMS DRIVEN BY FRACTIONAL BROWNIAN MOTIONS}

\author{ANHUI GU}

\address{Department of Mathematics, Shanghai Normal University,
Shanghai
200234, PR China \\ and College of Science, Guilin University of
Technology, \\
Guangxi,  Guilin 541004, PR China\\
gahui@glite.edu.cn}

\maketitle

\begin{history}
\received{January 8, 2012};
\revised{March 25, 2012}
\end{history}

\begin{abstract}
This paper is devoted to considering the stochastic lattice
dynamical systems (SLDS) driven by fractional Brownian motions with
Hurst parameter bigger than $1/2$. Under usual dissipativity
conditions these SLDS are shown to generate a random dynamical
system for which the existence and unique of a random attractor is
established. Furthermore, the random attractor is in fact a
singleton sets random attractor.
\end{abstract}


\begin{multicols}{2}
\section{Introduction}
\noindent The purpose of this paper is to investigate the
long-term behavior for the following SLDS
\begin{eqnarray} \label{1}
\frac{du_i(t)}{dt}&=&\kappa(u_{i-1}-2u_i+u_{i+1})-\lambda u_i
+f_i(u_i)\nonumber\\
&&+g_i+\sigma_i\frac{d\beta_i^H(t)}{dt}, ~ i\in \mathbb{Z},
\end{eqnarray}
where $\mathbb{Z}$ denotes the integer set,
$u=(u_i)_{i\in \mathbb{Z}}\in \ell^2 $,
$\kappa$ and $\lambda$ are positive constants, $f_i$ are smooth functions
satisfying some dissipative conditions,
$g=(g_i)_{i\in \mathbb{Z}}\in \ell^2$,
$\sigma=(\sigma_i)_{i\in \mathbb{Z}}\in \ell^2$,
and $\{\beta_i^H: i\in \mathbb{Z}\}$ are independent two-sided
fractional Brownian motions (fBms) with Hurst parameter $H\in (1/2, 1)$.

Recently, the dynamics of infinite lattice dynamical systems have
drawn much attention of mathematicians and physicists, see e.g.
\cite{BLW, BLL, Zhou1, Zhou2, Zhou3, ZS, Wang, LS1, LS2, Huang, WLW,
CL, ZZ1, WLX, HSZ, Han1, Han2} and the references therein. Since
most of the realistic systems involve noises which may play an
important role as intrinsic phenomena rather than just compensation
of defects in deterministic models, SLDS then arise naturally while
these random influences or uncertainties are taken into account.

Since Bates et al. \cite{BLL} initiated the study of SLDS,
many works have been done regarding the existence of global random
attractors for SLDS with white noises on lattices $\mathbb{Z}$
(see e.g. \cite{LS1, LS2, Huang, WLW, CL, ZZ1, WLX}).
Later, the existence of global random attractors have been extended to
other SLDS with additive white noises, for example, first-order
SLDS on $\mathbb{Z}^k$ \cite{LS1}, stochastic Ginzburg-Landau
lattice equations \cite{LS2}, stochastic FitzHugh-Nagumo lattice
equations \cite{Huang, WLW}, second-order stochastic lattice
systems \cite{WLX} and  the first (or second)-order SLDS with a
multiplicative white noise \cite{CL, Han1}. Zhao and Zhou \cite
{ZZ1} gave some sufficient conditions for the existence of a global
random attractor for general SLDS in the non-weighted space $\mathbb{R}$
of infinite sequences and provided an application to damped
sine-Gordon lattice system with additive noises. Very recently, Han,
et al. \cite{HSZ} provided some sufficient conditions for
the existence of global compact
random attractors for general random dynamical systems in weighted
space $\ell_\rho ^p$ $%
(p\geqslant 1)$ of infinite sequences, and their results are applied to
second-order SLDS in \cite{Han2}.

However, as can be seen that all the works above are considered in
the frameworks of the classical It$\ddot{\rm o}$ theory of Brownian
motion. There are no results on these systems when they are
perturbed by a fractional Brownian motion (fBm) to the best of our
knowledge. FBm appears naturally in the modeling of many complex
phenomena in applications when the systems are subject to ``rough"
external forcing. An fBm is a stochastic process which differs
significantly from the standard Brownian motion and
semi-martingales, and other classically used in the theory of
stochastic process. As a centered Gaussian process, it is
characterized by the stationarity of its increments and a medium- or
long-memory property. It also exhibits power scaling with exponent
$H$. Its paths are H$\ddot{\rm o}$lder continuous of any order
$H'\in (0, H)$. An fBm is not a semi-martingale nor a Markov
process. So an fBm is the good candidate to model random long term
influences in climate systems, hydrology, medicine and physical
phenomena. For more details on fBm, we can refer to the books
\cite{BHOZ, Mishura} for its further development.

The goal of this article is to establish the existence of a random
attractor for SLDS with the nonlinearity $f$ under some dissipative
conditions and driven by fBms with Hurst parameter $H\in (1/2, 1)$.
By borrowing the main ideas of the proofs in \cite{GKN}, we firstly
define a random dynamical system by using the possibility of a
pathwise interpretation of the stochastic integral with respect to
the fBms. This method is based on the fact that a stochastic
integral with respect to an fBm with Hurst parameter $H\in (1/2, 1)$
can be defined by a generalized pathwise Riemann-Stieltjes integral
(see e.g. \cite{Zahle, DU, NR, TTV}). And then we show the existence
of a pullback absorbing set for the random dynamical system which
achieved by means of a fractional Ornstein-Uhlenbeck transformation
and Gronwall lemma. Since every trajectory of the solutions of
system \eqref{1} cannot be differentiated, we have to consider the
difference between any two solutions among them, which is pathwise
differentiable (see \cite{GKN}). Due to the stationarity of the
fractional Ornstein-Uhlenbeck solution, we get the random attractor
finally. All solutions converge pathwise to each other, so the
random attractor is proven to be a singleton sets random attractor.

The paper is organized as follows. In Sec. 2, we recall some basic
concepts on random dynamical systems. In Sec.3, we give a unique
solution to system \eqref{1} and make sure that the solution
generates a random dynamical system. We establish the main result,
that is, the random attractor generated by equation \eqref{1} turns
out to be a singleton sets random attractor in Sec.4.

\section{Random dynamical systems and Random attractor}
In this section, we introduce some basic concepts related to random
dynamical systems and random attractor, which are taken from
\cite{CDF, Arnold}.

Let $(\mathbb{E}, \|\cdot\|_{\mathbb{E}})$ be a separable Hilbert space and
$(\Omega, \mathcal{F}, \mathbb{P})$ be a probability space.

\begin{definition} \label{MDS}
A metric dynamical system $(\Omega, \mathcal{F}, \mathbb{P},
\theta)$ with two-sided $\mathbb{R}$ consists of a measurable flow
\begin{equation*}
\theta: (\mathbb{R}\times \Omega, \mathcal{B}(\mathbb{R})\otimes \mathcal{F})
\rightarrow (\Omega, \mathcal{F}),
\end{equation*}
where the flow property for the mapping $\theta$ holds for the
partial mappings $\theta_t=\theta(t, \cdot)$:
\begin{equation*}
\theta_t\circ \theta_s=\theta_t\theta_s=\theta_{t+s}, \ \
\theta_0={\rm i d}_{\Omega}
\end{equation*}
for all $s, t\in \mathbb{R}$, and $\theta \mathbb{P}=\mathbb{P}$ for
all $t\in \mathbb{R}$.
\end{definition}

\begin{definition}
\label{RDS} A continuous random dynamical system
(RDS) $\varphi $ on $\mathbb{E}$ over $(\Omega ,\mathcal{F},
\mathbb{P},(\theta
_t)_{t\in \mathbb{R}})$ is a $(\mathcal{B}(\mathbb{R}^{+})\times \mathcal{F}%
\times \mathcal{B}(\mathbb{E}),\mathcal{B}(\mathbb{E}))$-measurable
mapping and satisfies

(i)  $\varphi (0,\omega )$ is the identity on $\mathbb{E}$;

(ii)  $\varphi (t+s,\omega )=\varphi (t,\theta _s\omega )\circ \varphi
(s,\omega )$ for all $s,$ $t\in \mathbb{R}^{+}$, $\omega \in \Omega $;

(iii)  $\varphi (t,\omega )$ is continuous on $\mathbb{E}$ for all
$(t,\omega )\in \mathbb{R}^{+}\times \Omega $.
\end{definition}

A universe $\mathcal{D}=\{D(\omega), \omega \in \Omega\}$ is a
collection of nonempty subsets $D(\omega)$ of $\mathbb{E}$
satisfying the following inclusion property: if $D\in \mathcal{D}$
and $D'(\omega)\subset D(\omega)$ for all $\omega \in \Omega$, then
$D'\in \mathcal{D}$.

\begin{definition}
A family of $\mathcal{A}=\{A(\omega), \omega \in \Omega\}$ of
nonempty measurable compact subsets $\mathcal{A}(\omega)$ of
$\mathbb{E}$ is called $\varphi$- invariant if $\varphi(t, \omega,
\mathcal{A}(\omega))=\mathcal{A}(\theta_t\omega)$ for all $t\in
\mathbb{R^+}$ and is called a random attractor if in addition it is
pathwise pullback attracting in the sense that
$$H_d^*(\varphi(t, \theta_{-t}\omega, D(\theta_{-t}\omega)),
\mathcal{A}(\omega))
\rightarrow 0 \ \ \mbox{as} \ \ t\rightarrow\infty$$
\end{definition}
for all $D\in \mathcal{D}$. Here $H_d^*$ is the Hausdorff
semi-distance on $\mathbb{E}$.

\begin{definition}
A random variable $r: \Omega\rightarrow \mathbb{R}$ is called tempered if
\begin{equation*}
\lim_{t\rightarrow\pm \infty}\frac{\log|r(\theta_t\omega)|}{|t|}=0
\ \ \mathbb{P}-a.s.
\end{equation*}
and a random set $\{D(\omega), \omega \in \Omega\}$ with
$D(\omega)\in \mathbb{E}$ is called tempered if it is contained in
the ball $\{x\in \mathbb{R}: |x|\leq r(\omega)\}$, where $r$ is a
tempered random variable.
\end{definition}

Here we will always work with the attracting universe given by the
tempered random sets.

\begin{definition}
A family $\hat{B}=\{B\mathcal{(\omega)}, \omega \in \Omega\}$ is
said to be pullback absorbing if for every $D(\omega)\in
\mathcal{D}$, there exists $T_{D}(\omega)\geq 0$ such that
\begin{equation}\label{2}
\varphi(t, \theta_{-t}\omega, D(\theta_{-t}\omega))\subset B(\omega) \ \
\forall t\geq T_{D}(\omega).
\end{equation}
\end{definition}

\begin{theorem}
\label{theorem 1}(See \cite{Schmalfuss, CDF, Arnold}.) Let $(\theta,
\varphi)$ be a continuous RDS on $\Omega\times \mathbb{E}$. If there
exists a pullback absorbing family $\hat{B}=\{B\mathcal{(\omega)},
\omega \in \Omega\}$ such that, for every $\omega \in \Omega$,
$B(\omega)$ is compact and $B(\omega)\in \mathcal{D}$, then the RDS
$(\theta, \varphi)$ has a random attractor

\[
\mathcal{A}(\omega )=\bigcap_{\tau>0}\overline{%
\bigcup_{t\geqslant \tau }\varphi (t,\theta _{-t}\omega )B(\theta
_{-t}\omega )}.
\]
\end{theorem}
Note that if the random attractor consists of singleton sets, i.e.
$\mathcal{A}(\omega)=\{u^*(\omega)\}$ for some random variable
$u^*$, then $u^*(t)(\omega)=u^*(t)(\theta_t\omega)$ is a stationary
stochastic process.

\section{SLDS with FBms}
We now firstly introduce as an example of metric dynamical system a
special noise that is called fractional Brownian motion. Given $H\in
(0, 1)$, a continuous centered Gaussian process $\beta^H(t), t\in
\mathbb{R}$, with the covariance function
\begin{equation*}
\mathbf{E} \beta^H(t)\beta^H(s)=\frac{1}{2}(|t|^{2H}+|s|^{2H}-|t-s|^{2H}),
\ \ t, s
\in \mathbb{R}
\end{equation*}
is called a two-sided one-dimensional fractional Brownian motion,
and $H$ is the Hurst parameter. For $H=1/2$, $\beta$ is a standard
Brownian motion, while for $H\neq 1/2$, it is neither a
semimartingale nor a Markov process. Moreover,
\begin{equation*}
\mathbf{E}|\beta^H(t)-\beta^H(s)|^2=|t-s|^{2H}, \ \ \mbox{for all}\ \
s, t\in \mathbb{R}.
\end{equation*}

Here, we assume that $H\in (1/2, 1)$ throughout the paper. When
$H\in (0, 1/2)$ we cannot define the stochastic integral by a
generalized Stieljes integral and, therefore, dealing with such
values of the Hurst parameter seems to be much more complicated. It
is worth mentioning that when $H=1/2$ the fBm becomes the standard
Wiener process, the random dynamical system generated by SLDS has
been studied in \cite{BLL}.

Using the definition of $\beta^H(t)$, Kolmogorov's theorem ensures
that $\beta^H$ has a continuous version, and almost all the paths
are H$\ddot{\rm o}$lder continuous of any order $H'\in (0, H)$ (see
\cite{Kunita}). Thus, we can consider the canonical interpretation
of an fBm: denote $\Omega=C_0(\mathbb{R}, \ell^2)$, the space of
continuous functions on $\mathbb{R}$ with values in $\ell^2$ such
that $\omega(0)=0$, equipped with the compact open topology. Let
$\mathcal{F}$ be the associated Borel-$\sigma$-algebra and
$\mathbb{P}$ the distribution of the fBm $\beta^H$, and
$\{\theta_t\}_{t\in \mathbb{R}}$ be the flow of Wiener shifts such
that
\begin{equation*}
\theta_t\omega(\cdot)=\omega(\cdot+t)-\omega(t), \ \ t\in \mathbb{R}.
\end{equation*}
Due to \cite{MS, GLS}, we know that the quadruple $(\Omega, \mathcal{F},
\mathbb{P}, \theta)$ is a
metric dynamical system which is ergodic. Furthermore, it holds that
\begin{eqnarray}\label{3}
&&\beta^H(\cdot, \omega)=\omega(\cdot), \nonumber\\
 \beta^H(\cdot, \theta_s\omega)&=&\beta^H(\cdot+s, \omega)
-\beta^H(s, \omega)\nonumber\\
&=&\omega(\cdot+s)-\omega(s).
\end{eqnarray}

We consider the SLDS
\begin{eqnarray} \label{4}
\frac{du_i(t)}{dt}&=&\kappa(u_{i-1}-2u_i+u_{i+1})-\lambda u_i+f_i(u_i)
+g_i\nonumber\\
&&+\sigma_i\frac{d\omega_i(t)}{dt}, ~ i\in \mathbb{Z},
\end{eqnarray}
where $\mathbb{Z}$ denotes the integer set, $u=(u_i)_{i\in
\mathbb{Z}}\in \ell^2$, $\kappa$ and $\lambda$ are positive
constants, $f_i$ are smooth functions satisfying some dissipative
conditions, $g=(g_i)_{i\in \mathbb{Z}}\in \ell^2$,
$\sigma=(\sigma_i)_{i\in \mathbb{Z}}\in \ell^2$, and
$\{\omega_i=\beta^H_i: i\in \mathbb{Z}\}$ are independent two-sided
fractional Brownian motions with Hurst parameter $H\in (1/2, 1)$.
Then, \eqref{4} can be understood as the pathwise Riemann-Stieltjes
integral equations
\begin{eqnarray}\label{5}
u_i(t)&=&u_i(0)+\int_0^t(\kappa(u_{i-1}(s)-2u_i(s)+u_{i+1}(s))
\nonumber\\
&&-\lambda u_i(s)+f_i(u_i(s))
+g_i)ds+\sigma_i\omega_i(t), ~ i\in \mathbb{Z}.\nonumber\\
\end{eqnarray}
Let $e^i\in \ell^2$ denote the element having 1 at position $i$ and
0 at other components. Then
\begin{equation}\label{6}
W(t)\equiv W(t, \omega)=\sum_{i\in \mathbb{Z}}\sigma_i\omega_i(t)e^i \ \
\mbox{with} \ \ (\sigma_i)_{i\in \mathbb{Z}}\in \ell^2,
\end{equation}
is the special noise with values in $\ell^2$ defined on the
probability space $(\Omega, \mathcal{F}, \mathbb{P})$.

{\bf Assumptions on the nonlinearity $f_i$.}  Assume $f_i:
\mathbb{R} \rightarrow \mathbb{R}$ to be continuously
differentiable. Let $f$ be the Nemytski operator associated with
$f_i$, i.e. for all $x=(x_i)_{i\in \mathbb{Z}}\in \ell^2$, let
$f(x)=(f_i(x_i))_{i\in \mathbb{Z}}$. Then we have $f(x)\in \ell^2$
(see \cite{BLL}). Assume $f$ satisfies a one-sided dissipative
Lipschitz condition
\begin{equation}\label{7}
\langle x-y, f(x)-f(y)\rangle \leq -L|x-y|^2, \ \ \mbox{for all} \ \
x, y\in \ell^2,
\end{equation}
where $L>0$, and a polynomial growth condition along with its
derivative, i.e.,
\begin{equation}\label{8}
|f(x)|+|Df(x)|\leq K(1+|x|^p), \ \ \mbox{for all} \ \ x\in \ell^2,
\end{equation}
$p\geq 1$ and a positive constant $K$. Here, we remark that we could
consider more general dissipativity conditions, which would lead to
nontrivial setvalued random attractors. However, to avoid rather
technical details,
 we will restrict here to the dissipativity
conditions \eqref{7} and \eqref{8}.

For $x=(x_i)_{i\in \mathbb{Z}}\in \ell^2$, define $A, B, B^*$ to be
linear operators from $\ell^2$ to $\ell^2$ as follows:
\begin{eqnarray*}(Ax)_i&=&-x_{i-1}+2x_i-x_{i+1}, \nonumber\\
(Bx)_i&=&x_{i+1}-x_i, \ \ (B^*x)_i=x_{i-1}-x_i,\ \ i\in \mathbb{Z}.
\end{eqnarray*}
It is easy to show that $A=BB^*=B^*B$, $\langle B^* x, x'\rangle
=\langle x, B x'\rangle$
for all $x, x'\in \ell^2$,
which implies that $\langle Ax, x\rangle\geq 0$.

The system \eqref{4} with initial values $u_0=(u_{0, i})_{i\in
\mathbb{Z}}\in \ell^2$ may be written as an equation in $\ell^2$
\begin{eqnarray}\label{9}
u(t)&=&u_0+\int_0^t(-\kappa Au(s)-\lambda u(s)+f(u(s))
+g)ds\nonumber\\&&
+W(t)\nonumber\\
&:=&u_0+\int_0^tG(u(s))ds+W(t),\nonumber\\
&& ~ t\geq 0, ~\omega\in \Omega,
\end{eqnarray}
where $G(u(t)):=-\kappa Au(t)-\lambda u(t)+f(u(t))
+g$.

\begin{proposition}\label{existence}
Let the above assumptions on $f$ be satisfied and $T>0$. Then system
\eqref{9} has a unique pathwise solution $u=(u(t))_{t\geq 0}$. In
addition that the solution satisfies
\begin{eqnarray} \label{10}
\sup_{t\in [0, T]}|u(t)|&\leq& M(1+|u_0|+\sum_{i\in \mathbb{Z}}
\sup_{\tau\in [0, T]}
|\beta_i^H(\tau)|\nonumber\\
&&+\sum_{i\in \mathbb{Z}}\sup_{\tau\in [0, T]}|\beta_i^H(\tau)|^{p}+|g|),
\end{eqnarray}
for all $T>0$, where $M$ is a positive constant and is independent
of $T$.
\end{proposition}

\begin{proof}
Let $v(t)=u(t)-W(t), t\geq 0$, then system \eqref{9} has a solution
$u=(u(t))_{t\geq 0}$ for all $\omega\in \Omega$ if and only if the
following equation
\begin{eqnarray}\label{11}
v(t)&=&u_0+\int_0^t(-\kappa Av(s)-\lambda v(s)+f(v(s)+W(s))
\nonumber\\
&&+g-\kappa AW(s)-\lambda W(s))ds.
\end{eqnarray}
has a unique pathwise solution for $t\in [0, T]$. However, since the
integrand is pathwise continuous, the fundamental theorem of
calculus says that the left hand side of \eqref{11} is pahtwise
differentiable. Thus, for fixed $\omega\in \Omega$, equation
\eqref{10} is the pathwise random ordinary differential equation
(RODE)
\begin{eqnarray}\label{12}
\frac{dv(t)}{dt}(\omega)&=&-\kappa Av(t)-\lambda v(t)+f(v(t)+W(t))
\nonumber\\
&&+g-\kappa AW(t)-\lambda W(t)\nonumber\\
&:=&\tilde{G}(t, v(t))(\omega):=\tilde{G}(v(t)+W(t)), \nonumber\\
&&~~~~\mbox{for}\ \  t\geq 0, \ \  v_0(\omega)=u_0(\omega).
\end{eqnarray}
Since $\tilde{G}: [0, \infty)\times \ell^2\rightarrow \ell^2$ is
continuous in $t$ and continuous differentiable in $v$, this RODE
has a unique local solution in a small interval $[0, \tau(\omega)]$,
which means that \eqref{9} has a unique local solution in the same
small interval $[0, \tau(\omega)]$, see e.g. Theorem 2.1.4 in
\cite{SH}.

To see that \eqref{9} has a unique solution for every $t\geq 0$, we
prove at first the a priori estimate \eqref{10}. Suppose that
$u=(u(t))_{t\geq 0}$ solves \eqref{9} on the interval $[0, T]$. This
implies that the process $v(t)=u(t)-W(t)$ solves \eqref{11} on the
interval $[0, T]$. By the one-sided dissipative Lipschitz condition
and H$\ddot{\rm o}$lder inequality, we obtain
\begin{eqnarray*}
&&2|v(t)|\frac{d|v(t)|}{dt}=\frac{d}{dt}|v(t)|^2\nonumber\\
&=&2\langle v(t), \tilde{G}(v(t)+W(t))\rangle\nonumber\\
&=&2\langle v(t), -\kappa Av(t)-\lambda v(t)+f(v(t)+W(t))
\nonumber\\
&&+g-\kappa AW(t)-\lambda W(t)\rangle\nonumber\\
&\leq& -2(\lambda+L)|v(t)|^2+2|f(W(t))||v(t)|\nonumber\\
&&+2(|g|-|\kappa AW(t)|-|\lambda W(t)|)|v(t)|,
\end{eqnarray*}
that is
\begin{eqnarray*}
\frac{d|v(t)|}{dt}
&\leq& -(\lambda+L)|v(t)|+c_0(|g|+\sup_{\tau\in [0, T]}|W(\tau)|
\nonumber\\
&&+\sup_{\tau\in [0, T]}|W(\tau)|^{p}),
\ \ t\in [0, T],
\end{eqnarray*}
where $c_0$ is a positive constant depends on $\kappa, \lambda, K,
\max_{j\in \mathbb{Z}}|\sigma_j|$
and $\max_{j\in \mathbb{Z}}|g_j|$.
By Gronwall lemma it yields that
\begin{eqnarray*}
|v(t)|&\leq& |u_0|e^{-\lambda t}+\frac{c_0}{\lambda}(1-e^{-\lambda t})(|g|
+\sup_{\tau\in [0, T]}|W(\tau)|\nonumber\\
&&+\sup_{\tau\in [0, T]}|W(\tau)|^{p}),
\ \ t\in [0, T].
\end{eqnarray*}
Since $\lambda>0$, we have
\begin{eqnarray*}
|v(t)|&\leq& |u_0|+\frac{c_0}{\lambda}(|g|+\sup_{\tau\in [0, T]}|W(\tau)|
\nonumber\\
&&+\sup_{\tau\in [0, T]}|W(\tau)|^{p}),
\end{eqnarray*}
that is
\begin{eqnarray*}
\sup_{t\in [0, T]}|v(t)|&\leq& M(|g|+\sup_{\tau\in [0, T]}|W(\tau)|
\nonumber\\
&&+\sup_{\tau\in [0, T]}|W(\tau)|^{p}),
\ \ t\in [0, T],
\end{eqnarray*}
where the constant $M$ depends on $\kappa, \lambda, K, \max_{j\in
\mathbb{Z}}|\sigma_j|$ and $\max_{j\in \mathbb{Z}}|g_j|$. The result
of \eqref{10} follows then by using the relation $v(t)=u(t)+W(t)$
for $t\in [0, T]$. As a consequence of estimate \eqref{10} the local
unique solution to \eqref{9} can be extended to a global unique
solution (see e.g. Theorem 2.1.4 in \cite{SH}).
\end{proof}

We consider \eqref{1} with the linear
drift term $G(u(t))=-\lambda u(t)$, that is
\begin{equation}\label{16}
du(t)=-\lambda u(t)dt+dW(t), \ \ \lambda>0,
\end{equation}
which is called the fractional Ornstein-Uhlenbeck process, where
$W(t)$ denotes a one-dimensional fractional Brownian motion defined
in \eqref{6}. It has the explicit solution
\begin{equation}\label{17}
u(t)=u_0e^{-\lambda t}+e^{-\lambda t}\int_0^te^{\lambda s}dW(s).
\end{equation}
Taking the pathwise pullback limit, we get the stochastic stationary
solution
\begin{equation}\label{18}
\bar{u}(t)=e^{-\lambda t}\int_{-\infty}^te^{\lambda s}dW(s), \ \
t\in \mathbb{R},
\end{equation}
which is called the fractional Ornstein-Uhlenbeck solution. Based on
Lemma 2.6 in \cite{MS} and Lemma 1 in \cite{GKN}, we have the
property
\begin{lemma} \label{property}
For all $\omega\in \Omega$ the Riemann-Stieltjes integrals
\begin{equation*}
e^{-\lambda t}\int_{-\infty}^te^{\lambda s}dW(s), \ \ t\in \mathbb{R},
\end{equation*}
are well defined in $\ell^2$. Moreover, for all $\omega\in \Omega$,
we have
\begin{equation*}
|e^{-\lambda t}\int_{-\infty}^te^{\lambda s}dW(s)|\leq 4\rho(\omega)(1+|t|)^2, \ \
t\in \mathbb{R},
\end{equation*}
where the random constant $\rho(\omega)>0$ and satisfying
\begin{equation}\label{19}
|W(t)|\leq \rho(\omega)(1+|t|^2), \ \ \mbox{for} \ \  t\in \mathbb{R},
\ \ \omega\in \bar{\Omega}.
\end{equation}
Here, $\bar{\Omega}\in \mathcal{F}$ is a $(\theta_t)_{t\in
\mathbb{R}}$-invariant set of full measure.
\end{lemma}

\begin{proof} We know that $\lambda>0$, $(\sigma_i)_{i\in \mathbb{Z}}\in \ell^2$.
The proof is similar to Lemma 2.6 in \cite{MS} and
Lemma 1 in \cite{GKN}, thus we omit it here.

\end{proof}

Since the mapping of $\theta$ on $\bar{\Omega}$ has the same
properties as the original one if we choose the trace
$\sigma$-algebra with respect to $\bar{\Omega}$ to be denoted also
by $\mathcal{F}$, we can change our metric dynamical system with
respect to $\bar{\Omega}$, and still denoted by the symbols
$(\Omega, \mathcal{F}, \mathbb{P},(\theta_{t})_{t\in \mathbf{R}})$.

Now, we will verify that the solution of system \eqref{1} generates
a continuous RDS.

\begin{proposition} \label{CRDS}
The solution of \eqref{1} determinants a continuous random dynamical
system $\varphi: \mathbb{R^+}\times \Omega\times \ell^2 \rightarrow
\ell^2$, which is given by
\begin{eqnarray}\label{13}
\varphi(t, \omega, u_0)&=&u_0+\int_0^t (-\kappa Au(s)-\lambda u(s)
\nonumber\\
&&+f(u(s))
+g)ds+W(t, \omega).
\end{eqnarray}
\end{proposition}

\begin{proof}
Note that \eqref{3} is satisfied for $\omega\in \Omega$ and from the
definition of $(\theta_t)_{t\in \mathbb{R}}$, we have the property
\begin{equation}\label{14}
W(\tau+t, \omega)=W(\tau, \theta_{t}\omega)+W(t, \omega) \ \ \mbox{for all}
\ \ t, \tau\in \mathbb{R}.
\end{equation}
From Proposition \ref{existence} we know that $\varphi$ solves
\eqref{1}, thus $\varphi$ is measurable and satisfies $\varphi(0,
\omega, \cdot)=\rm id_{\ell^2}$. It remains to verify the cocycle
property in Definition \ref{RDS}. Let $t, \tau\in \mathbb{R^+},
\omega\in \Omega$ and $u_0\in \ell^2$, denote
$G(u(t))(\omega):=-\kappa Au(t)-\lambda u(t)+f(u(t)) +g$, it yields
from \eqref{3} that
\begin{eqnarray} \label{15}
&&\varphi(t+\tau, \omega, u_0)\nonumber\\
&=&u_0+\int_0^{t+\tau} G(u(s))(\omega)ds+W(t+\tau, \omega)\nonumber\\
&=& u_0+\int_0^t G(u(s))(\omega)ds+W(t, \omega)\nonumber\\
&&+\int_t^{t+\tau} G(u(s))(\omega)ds+W(\tau, \theta_t\omega)\nonumber\\
&=& u(t)+\int_0^{\tau} G(u(s))(\theta_t\omega)ds
+W(\tau, \theta_t\omega)\nonumber\\
&=& \varphi(\tau, \theta_{t}\omega, \cdot)\circ \varphi(t, \omega, u_0),
\end{eqnarray}
which completes the proof.
\end{proof}

\section{Existence of a Random attractor}

In this section, we will prove the existence of a random attractor
for the RDS defined in Proposition \ref{CRDS}. The main proof is
based on the first part in section 4 in \cite{GKN}.

Let $u, w$ be any two solutions of system \eqref{1}. Their sample
paths are not differentiable, but the difference satisfies pathwise
\begin{eqnarray*}
u(t)-w(t)&=&u_0-w_0+\int_0^t(-\kappa A(u(s)-w(s))\nonumber\\
&&-\lambda(u(s)-w(s))\nonumber\\
&&+(f(u(s))-f(w(s))))ds, \ \ t\geq 0,
\end{eqnarray*}
and again, since the integrand is pathwise continuous, the
fundamental theorem of calculus indicates that the left hand side is
pathwise differentiable and satisfies
\begin{eqnarray}\label{20}
&&\frac{d}{dt}(u(t)-w(t))\nonumber\\
&=&-\kappa A(u(s)-w(s))-\lambda(u(s)-w(s))\nonumber\\
&&+(f(u(s))-f(w(s))), \ \ t\geq 0.
\end{eqnarray}
Now we use condition \eqref{7} again, we obtain from \eqref{20} that
\begin{eqnarray*}
&&\frac{d}{dt}|u(t)-w(t)|^2\nonumber\\
&=&2\langle u(t)-w(t), -\kappa A(u(s)-w(s))\nonumber\\
&&-\lambda(u(s)-w(s))+(f(u(s))-f(w(s)))\rangle\nonumber\\
&\leq& -2\lambda |u(t)-w(t)|^2.
\end{eqnarray*}
Thus pathwise we have
\begin{equation*}
|u(t)-w(t)|\leq |u_0-w_0|e^{-\lambda t}\rightarrow 0, \ \
\mbox{as} \ \ t\rightarrow\infty.
\end{equation*}
That is to say that all solutions converge pathwise forward to each
other in time.

Now, we consider the difference $u(t)-\bar{u}(t)$. This is also
pathwise differentiable, since the paths are continuous and satisfy
the integral equation
\begin{eqnarray*}
u(t)-\bar{u}(t)&=&u_0-\bar{u}_0+\int_0^t(-\kappa Au(s)
-\lambda(u(s)\nonumber\\
&&-\bar{u}(s))+f(u(s)))ds,
\end{eqnarray*}
which is equivalent to the pathwise differential equation
\begin{eqnarray*}
&&\frac{d}{dt}(u(t)-\bar{u}(t))\nonumber\\
&=&-\kappa Au(s)-\lambda(u(s)-\bar{u}(s))+f(u(s)), t\geq 0.
\end{eqnarray*}
By using \eqref{7}, it follows that
\begin{eqnarray*}
\frac{d}{dt}|u(t)-\bar{u}(t)|^2&=&2\langle u(t)-\bar{u}(t),
-\kappa Au(s)\nonumber\\
&&-\lambda(u(s)-\bar{u}(s))+f(u(s))\rangle\nonumber\\
&\leq& -2\lambda |u(t)-\bar{u}(t)|^2\nonumber\\
&&+2|u(t)-\bar{u}(t)||f(\bar{u}(t))|,
\end{eqnarray*}
that is
\begin{eqnarray*}
\frac{d}{dt}|u(t)-\bar{u}(t)|&\leq& -\lambda |u(t)-\bar{u}(t)|
+|f(\bar{u}(t))|,
\end{eqnarray*}
and hence
\begin{eqnarray}\label{21}
&&|u(\omega)-\bar{u}(\omega)|\nonumber\\
&\leq& |u_0(\omega)-\bar{u}_0(\omega)|e^{-\lambda t}\nonumber\\
&&+e^{-\lambda t}\int_0^te^{\lambda s}|f(\bar{u}(s))(\omega)|ds.
\end{eqnarray}
Let us check that the family of balls centered on
$\bar{u}_0(\omega)$ with the random radius
\begin{equation}\label{22}
\varrho(\omega):=1+\int_{-\infty}^0e^{\lambda s}|f(\bar{u}(s))(\omega)|ds
\end{equation}
is a pullback absorbing family for the random dynamical system
generated by system \eqref{1}.

Due to the assumptions on $f$ and Lemma \ref{property}, the radius
defined in \eqref{22} is well defined. Now, by replacing $\omega$ by
$\theta_{-t}\omega$ in \eqref{21}, we get
\begin{eqnarray}\label{23}
&&|u(\theta_{-t}\omega)-\bar{u}(\theta_{-t}\omega)|\nonumber\\
&\leq& |u_0(\theta_{-t}\omega)-\bar{u}_0(\theta_{-t}\omega)|
e^{-\lambda t}\nonumber\\
&&+\int_0^te^{\lambda (s-t)}|f(\bar{u}(s))(\theta_{-t}\omega)|ds
\nonumber\\
&=& |u_0(\theta_{-t}\omega)-\bar{u}_0(\theta_{-t}\omega)|
e^{-\lambda t}\nonumber\\
&&+\int_{-t}^0e^{\lambda \tau}|f(\bar{u}(\tau))(\omega)|d\tau,
\ \ t\geq 0.
\end{eqnarray}
The last term in \eqref{23} due to $\bar{u}(s)(\theta_{-t}\omega)=
\bar{u}_0(\theta_{s-t}\omega)=\bar{u}(s-t)(\omega)$ which deduced
from that $(\bar{u}(t))_{t\in \mathbb{R}}$ is a stationary process.
The conclusion now follows for $t\rightarrow\infty$.

Because of the stationarity and Lemma \ref{property}, we have
$e^{-\lambda t}|\bar{u}_0(\theta_{-t}\omega)|
=e^{-\lambda t}|\bar{u}(-t)(\omega)|
\rightarrow 0$ as $t\rightarrow\infty$. Then
we have the pullback absorption
\begin{equation}\label{24}
|u(\theta_{-t}\omega)|\leq |\bar{u}_0(\omega)|
+\varrho(\omega), \ \ \forall t\geq T_{\mathcal{D}(\omega)}.
\end{equation}
As a consequence of Theorem \ref{theorem 1}, system \eqref{1} has a
random attractor $\mathcal{A}=\{A(\omega), \omega\in \Omega\}$. We
have know that all solutions converge pathwise to each other, so the
random attractor sets are in fact singleton sets
$\mathcal{A}=\{\tilde{u}_0(\omega)\}$, i.e. the random attractor is
formed by a stationary random process
$\tilde{u}(t)(\omega):=\tilde{u}_0(\theta_t\omega)$, which pathwise
attracts all other solutions in both forward and pullback senses.

Finally, we are in the position to state the result of existence of
a random attractor.
\begin{theorem}
Assume that the conditions on $f$ be satisfied. Then the random
dynamical system $\varphi$ has a unique random attractor.
Furthermore, the random attractor is in fact a singleton sets random
attractor.
\end{theorem}

\begin{remark}
Sometimes, for the need of demonstrating the relations between
$\frac{dv(t)}{dt}(\cdot)$ (or $G(u(t))(\cdot)$, $\tilde{G}(t,
v(t))(\cdot)$) and $\omega$ more explicitly, we will write
$\frac{dv(t)}{dt}(\omega)$ (or $G(u(t))(\omega)$, $\tilde{G}(t,
v(t))(\omega)$) instead if necessary.
\end{remark}

\begin{remark}
There are several differences between these two kinds of noises:
white noise and fBm for the existence of a random attractor of SLDS.
On one hand, in this work, the SLDS perturbed by additive fBms are
considered, which allow to transform the SLDS into random systems
that can be dealt with in a pathwise way. But for the multiplicative
noise, with no transformation into a random system, the existence of
random attractor for stochastic differential equations was consider
in \cite{GMS} based on considering a suitable sequence of stopping
times for the fractional Brownian motion. This method is different
from the traditional conjugacy method, which transforms a stochastic
equation into a differential equation with random coefficients but
without white noise. On the other hand, when SLDS perturbed by
either additive white noises (see e.g. \cite{BLL}) or a
multiplicative white noise (see e.g. \cite{CL})), the existing
random attractor is a (nontrivial) compact set of tempered random
bounded set, but the random attractor turns out to be a single sets
random attractor when the systems disturbed by additive fBms.
\end{remark}

\nonumsection{Acknowledgments} \noindent The author would like to
thank the anonymous referees for very helpful comments and
suggestions which largely improve the quality of the original
manuscript. This work has been partially supported by NSF of China
grants 11071165 and Guangxi Provincial Department of Research
Project grants 201010LX166.

\end{multicols}
\end{document}


\baselineskip 20pt

\begin{center}

\textbf{\Large Erratum: Random Attractors of Stochastic Lattice
Dynamical Systems Driven by Fractional Brownian Motions. [Int. J.
Bifurcation Chaos 23, 1350041 (2013)]}

\vskip 0.5cm

{\large  Anhui Gu}

\vskip 0.3cm

\end{center}

\vspace{10pt}

\baselineskip 18pt
\begin{multicols}{2}
In the recent paper \cite{Gu} we have obtained the existence of a
random attractor which turns out to be a singleton sets random
attractor. The result relies on the theory in \cite{Arnold}, that is
the existence of a compact pullback absorbing set, which just deals
with random dynamical systems in finite dimensional cases.

However, in the first-order lattice dynamical systems driven by
fractional Brownian motions, we considered it in the framework of
infinite dimensional square-summable sequences space $\ell^2$. The
main result need few changes because we cannot obtain it based on
Theorem 1 in \cite{Gu} but on the assumptions of the nonlinear
function $f$. The random dynamical system generated by the system
has a unique random equilibrium, all solutions converge pathwise to
each other, so the random attractor, which consists of a unique
random equilibrium, is proven to be a singleton sets random
attractor. With the misused theory several important parts of the
paper should be corrected as:

\begin{itemize}
\item{} The last sentence in the Abstract section should be ``In our case,
the random dynamical system has a unique random equilibrium, which
constitutes a singleton sets random attractor."

\item{} In page 7, from line 15 to 29, the main statement should
be adjusted to ``Assume that the conditions on $f$ are satisfied.
Then the random dynamical system $\varphi$ has a unique random
equilibrium, which constitutes a singleton sets random attractor."
\end{itemize}

\end{multicols}{2}